\documentclass[a4paper]{amsart}

\usepackage{amssymb}

\newcommand{\intt}{\mathop{\mathrm{int}}\nolimits}
\newcommand{\im}{\mathop{\mathrm{Im}}\nolimits}
\newcommand{\re}{\mathop{\mathrm{Re}}\nolimits}
\newcommand{\e}{\varepsilon}
\newcommand{\DD}{\hbox{D\kern-.73em\raise.25ex\hbox{-}\raise-.25ex\hbox{ }}}
\newcommand{\dD}{\hbox{d\kern-.33em\raise.75ex\hbox{-}\raise-.25ex\hbox{}}}
\newenvironment{pf}{\begin{proof}}{\end{proof}}

\usepackage{amsmath}
\usepackage{amsfonts}
\usepackage{amssymb}

\newtheorem{Theorem}{Theorem}[section]
\newtheorem{Lemma}[Theorem]{Lemma}
\newtheorem{Proposition}[Theorem]{Proposition}
\newtheorem{Corollary}[Theorem]{Corollary}

\theoremstyle{definition}
\newtheorem{Definition}[Theorem]{Definition}

\theoremstyle{remark}
\newtheorem{Remark}[Theorem]{Remark}

\begin{document}

\date{}


\title[TVS-cone metric spaces as a special case of metric spaces]{TVS-cone metric spaces as a special case of metric spaces}

\author{Ivan D.~Aran\dD elovi\'c}

\address{University of Belgrade - Faculty of Mechanical
Engineering, Kraljice Marije 16, 11000 Beograd, Serbia}

\email{iarandjelovic@mas.bg.ac.rs}

\author{Dragoljub J.~Ke\v cki\'c}

\address{University of Belgrade - Faculty of Mathematics, Studentski trg 16, 11000 Beograd, Serbia}

\email{keckic@matf.bg.ac.rs}

\thanks{The first author was partially supported by MNZZS Grant, No. 174002, Serbia.
The second author was partially supported by MNZZS Grant, No. 174034, Serbia.}

\begin{abstract}

There have been a number of generalizations of fixed point results to the so called TVS-cone metric spaces,
based on a distance function that takes values in some cone with nonempty interior (solid cone) in some
topological vector space. In this paper we prove that the TVS-cone metric space can be equipped with a family
of mutually equivalent (usual) metrics such that the convergence (resp.\ property of being Cauchy sequence,
contractivity condition) in TVS sense is equivalent to convergence  (resp.\ property of being Cauchy
sequence, contractivity condition) in all of these metrics. As a consequence, we prove that if a topological
vector space $E$ and a solid cone $P$ are given, then the category of TVS-cone metric spaces is a proper
subcategory of metric spaces with a family of mutually equivalent metrics (Corollary \ref{categoty}). Hence,
generalization of a result from metric spaces to TVS-cone metric spaces is meaningless. This, also, leads to
a formal deriving of fixed point results from metric spaces to TVS-cone metric spaces and makes some earlier
results vague. We also give a new common fixed point result in (usual) metric spaces context, and show that
it can be reformulated to TVS-cone metric spaces context very easy, despite of the fact that formal
(syntactic) generalization is impossible. Apart of main results, we prove that the existence of a solid cone
ensures that the initial topology is Hausdorff, as well as it admits a plenty of convex open sets. In fact
such topology is stronger then some norm topology.
\end{abstract}



\maketitle


\section{Introduction}
\label{1}

Investigation of cone metric spaces (also known as $K$-metric space) was introduced by several Russian
authors in the middle of 20th century. They differ from usual metric spaces in the fact that the values of
distance functions are not positive real numbers, but elements of a cone in some normed space. L.-G.~Huang
and X.~Zhang \cite{HZ}, reintroduced such spaces, and also went further, defining convergent and Cauchy
sequences in the terms of interior points of the underlying cone. They proved some fixed point theorems for
this class of spaces. Du \cite{DU} generalizes this notion considering a cone in some locally convex
topological vector space $E$ instead of a normed space. Further generalization is given by I.~Beg, A.~Azam
and M.~Arshad, \cite{BAA} assuming that the distance function takes values in Hausdorff (not necessarily
locally convex) topological vector space. Such spaces are known as TVS-cone metric spaces.

Given a TVS-cone metric space $(X,d)$ (over topological vector space $E$), we shall prove that there exists a
family of metrics $\{d_\alpha\}$ (in the usual sense - with real values) that defines exactly the same
convergent and Cauchy sequences and such that inequality in TVS metric is equivalent to the same inequality
in all metrics $d_\alpha$. More precisely, the category of TVS-cone metric spaces is a proper subcategory of
metric spaces equipped with a family of mutually equivalent metrics. Hence, to say that some result is a
generalization of a result from metric spaces to TVS-cone metric spaces is same as saying that a result is
generalization from general topology to metric spaces, i.e, it is a nonsense. Therefore, a large number of
generalizations of fixed point results to TVS-cone metric spaces (or to its special case - cone metric
spaces), published in last several years are not really generalizations, but just a complicated way to
formulate a result that is a special case of an old one.

It has to be mentioned that this result was partially proved in \cite{KP}, \cite{KRR1} and \cite{AR}.

The main tool we used is the scalarization function, introduced by Du \cite{DU}. For any scalarization
function, we construct its induced norm on the underlying topological vector space $E$. Further, we show that
all such obtained norms are mutually equivalent, and that convergence in any of them is equivalent to TVS
convergence. The same is done for Cauchy sequences, continuous and uniformly continuous functions.

Also, we prove that these equivalent norms are continuous with respect to the initial topology, and hence,
they generate topology weaker then initial. As a consequence we derived that the existence of a solid cone in
topological vector space ensures that the topology is Hausdorff as well as it contains many convex open
neighborhood of zero.

In the last section, we give Theorem that allows automatic extension of fixed point results by formal
syntactic algorithm, and also illustrate, by Example, that the extension is not less easy if formal methods
are unavailable. This example includes new fixed point result on usual metric spaces.

We give the proofs of many basic statements initially proved in some other paper to make this note more
readable.

\section{Scalarization function and induced norms}
\label{2}

\begin{Definition} Let $E$ be a topological vector space (over the field $K=\mathbf R$ or $\mathbf C$),
and let $P\subseteq E$ be a closed set with following properties:
\begin{enumerate}

\item $P\neq\emptyset,\{0\}$;

\item $x,y\in P$, and $\lambda\ge0$ implies $x+y,\lambda x\in P$;

\item\label{d12} $P\cap -P=\{0\}$.

\end{enumerate}
Then $P$ is called a {\em cone}. If, in addition, its interior is nonempty, we say that $P$ is a {\em solid
cone}.
\end{Definition}

\begin{Definition}\label{poredak}
We define a partial ordering with respect to cone $P$. We write $x\le y$, and say "$x$ is less then $y$", if
$y-x\in P$, and $x\ll y$ ("$x$ is strongly less then $y$") for $y-x\in\intt P$. The latter, of course, if the
cone is solid. It is easy to verify that this relation is reflexive ($0\in P$), antisymmetric ($P\cap
-P=\{0\}$) and transitive (follows from subadditivity).

If ambiguity is possible we shall emphasize notations and use $\le_P$ and $\ll_P$.

In further we shall write tvs for topological vector space. The pair $(E,P)$ consisting of tvs $E$ and a
solid cone $P$ is called {\em ordered tvs}.
\end{Definition}

\begin{Remark} Sometimes, the third condition is omitted in the definition of a cone. In this case, cones
which satisfies \ref{d12}, are called {\em pointed cones} or {\em proper cones}.
\end{Remark}

The following Lemma establishes basic properties of a solid cone. The second one is the most important and
due to it, TVS-cone convergence does not make anything new.

\begin{Lemma}\label{basic} Let $(E,P)$ be an ordered tvs. Then:

1) If $x\in P$ and $y\in\intt P$ then $x+y\in\intt P$. Consequently, if $x\le y$ and $y\ll z$ (or $x\ll y$
and $y\le z$) then $x\ll z$;

2) For any $c,e\in\intt P$ there are $\delta,\eta>0$ such that $\delta c\ll e$ and $\eta e\ll c$,
\end{Lemma}

\begin{pf} 1) Since $y\in\intt P$ there is an open neighborhood $W$ of zero such that $y+W\subseteq P$. On
the other hand $P$ is closed with respect to addition and hence $x+y+W\subseteq P$ which means $x+y\in\intt
P$;

2) The set $\{x\in E~|~x\ll e\}$ is equal to $e-\intt P$; hence it is open. Also it contains the origin;
therefore it contains some neighborhood of zero $W$. Since $c/n\to0$, there is a positive integer $n$ such
that $c/n\in W$; Take $\delta=1/n$ to obtain $\delta c\ll e$.
\end{pf}

\begin{Definition} Let $e\in\intt P$. The nonlinear scalarization function $\xi_e:E\to E$ is defined by
$$\xi_e(y)=\inf\{t\in\mathbf{R}~|~y\in te-P\}=\inf\{t\in\mathbf R~|~y\le te\}=\inf M_{e,y}.$$
\end{Definition}

To see that $\xi_e$ is well defined we must check that the set $M_{e,y}$ is nonempty and bounded below.
Indeed, $e\in\intt P$, as well as $-e\notin P$, and both $\intt P$, $E\setminus P$ are open. Therefore, there
is a neighborhood $W$ of zero such that $e+W\subseteq P$, and $(-e+W)\cap P=\emptyset$. Obviously $-y/n\to
0$, as $n\to\infty$, and therefore, there is a positive integer $n$ such that $-y/n\in W$ and then $ne-y\in
P$ and $-ne-y\notin P$; hence $n\in M_{e,y}$, $-n\notin M_{e,y}$.

\begin{Remark}\label{rem1} In fact, since $P$ is closed, the set $M_{e,y}$ is also closed, and hence
\begin{equation}\label{skupM}
M_{e,y}=[\xi_e(y),+\infty).
\end{equation}
\end{Remark}

The following Lemma that establishes the basic properties of $\xi_e$ was proved (except statement 6)) in
\cite{DU}.

\begin{Lemma}\label{sf} For any $e\in\intt P$, the function $\xi_e$ has the following properties:

1) $\xi_e(0)=0$;

2) $y\in P$ implies $\xi_e(y)\ge 0$;

3) if $y_2\le y_1$, then $\xi_e(y_2)\le \xi_e(y_1)$ for any $y_1,y_2\in E$;

4) $\xi_e$ is subadditive on $E$.

5) $\xi_e$ is positively homogeneous on $E$.

6) $\xi_e$ is continuous on $E$;
\end{Lemma}

\begin{pf} 1) $M_{e,0}=\{t~|~te\in P\}=[0,+\infty)$;

2) Let $y\in P$. If $t<0$, and $te-y\in P$ then $te=te-y+y\in P$ which contradicts $te\in -P$;

3) If $y_2\le y_1$, then $y_1-y_2\in P$, and for each $t\in M_{e,y_1}$ it holds $y_1\le te$. Hence $y_2\le
te$;

4) From $x\in t_1e-P$ and $y\in t_2e-P$ it follows $$x+y\in t_1e+t_2e-(P+P)=(t_1+t_2)e-P.$$ So $\xi_e(x)\le
t_1$ and $\xi_e(y)\le t_2$ implies $\xi_e(x+y)\le t_1+t_2$. Hence $\xi_e(x+y)\le \xi_e(x)+\xi_e(y)$;

5) It is enough to check that $M_{e,\lambda y}=\lambda M_{e,y}$ for any $\lambda>0$.

6) First, we have $\xi_e^{-1}(\alpha,+\infty)=\{x~|~\xi_e(x)>\alpha\}=\{x~|~x\not\in \alpha e-P\}=\alpha
e-E\setminus P$, which is an open set, since $P$ is closed.

Next, let $x\in\xi_a^{-1}(-\infty,\alpha)$. Then $\xi_e(x)=\alpha-\delta$, for some $\delta>0$. The mapping
$y\mapsto\delta y/2$ is continuous, implying $\delta e/2\in\intt P$. Then there exists a balanced
neighborhood of zero $W$, such that $\delta e/2+W\subseteq\intt P$. For $y\in W$ we have $-y\in W$, as well,
and therefore $\delta e/2-y\in P$, implying $\xi_e(y)\le\delta/2$, and also, by subadditivity
$$\xi_e(x+y)\le\xi_e(x)+\xi_e(y)\le\alpha-\delta+\delta/2<\alpha.$$
Hence, $x+W\subseteq\xi_e^{-1}(-\infty,\alpha)$ which finishes the proof.
\end{pf}

\begin{Lemma}\label{l1} Let $(E,P)$ be an ordered tvs, let $e\in\intt P$, and let $\xi_e$ be its scalarization
function.

a) If $0\ll x\ll\lambda e$, for some real $\lambda>0$, then $0\le \xi_e(x),-\xi_e(-x)<\lambda$.

b) The function $E\ni x\mapsto \|x\|_e=\max\{|\xi_e(x)|,|\xi_e(-x)|\}$ is a norm on $E$ considered over real
field.

c) If $x\in P$ then $\|x\|_e=\xi_e(x)\ge0$;
\end{Lemma}

\begin{Remark} If $E$ is a linear space over $\mathbf C$, we can regard it as a linear space over $\mathbf R$.
Thus, $\|\cdot\|_e$ must satisfy $\|\lambda x\|=|\lambda|\,\|x\|$ only for real $\lambda$. In the rest of the
paper, we deal only with real norms.
\end{Remark}

\begin{pf}  a) The left hand side of the required inequality follows from Lemma \ref{sf} - 2). From Lemma
\ref{sf} - 3) it holds $\xi_e(x)\le\xi_e(\lambda e)=\lambda$. Since $\lambda e-x\in\intt P$, there exists a
neighborhood $W$ of $0$ such that $\lambda e-x+W\subseteq P$. Let $\delta>0$ be such that $-\delta e\in W$.
Then we have
\begin{equation}\label{prva}P\ni \lambda e-x-\delta e=(\lambda-\delta)e-x,
\end{equation}
and therefore $\xi_e(x)\le\lambda-\delta<\lambda$.

Further, from (\ref{prva}), we conclude that $(\delta-\lambda)e+x\notin P$, and hence $\xi_e(-x)>-\lambda$.

b) Obviously, $\|x\|_e\ge 0$, $\|0\|_e=0$ and $\|-x\|_e=\|x\|_e$. Subadditivity follows from subadditivity of
the functions $x\mapsto\xi_e(x)$, $x\mapsto \xi_e(-x)$, $|.|$ and $\max$, and the fact that at least one of
numbers $\xi_e(x+y)$, $\xi_e(-x-y)$ must be positive (the latter follows from
$0=\xi_e(x+y+(-x-y))\le\xi_e(x+y)+\xi_e(-x-y)$).

Since $\xi_e$ is positively homogeneous, for $\lambda>0$ we have $\xi_e(\lambda x)=\lambda\xi_e(x)$ and
$\xi_e(-\lambda x)=\lambda\xi_e(-x)$, implying $|\xi_e(\lambda x)|=\lambda|\xi_e(x)|$ and $|\xi_e(-\lambda
x)|=\lambda|\xi_e(-x)|$; hence $\|\lambda x\|_e=\lambda\|x\|_e$. If $\lambda<0$ then we use the previous case
$\lambda>0$ and $\|-x\|_e=\|x\|_e$ to obtain $\|\lambda
x\|_e=\|(-\lambda)(-x)\|_e=(-\lambda)\|-x\|_e=|\lambda|\,\|x\|_e$.

Finally, if $\|x\|_e=0$ then both $\xi_e(x)$ and $\xi_e(-x)$ are equal to zero. It follows from (\ref{skupM})
that $x,-x\in P$ implying $x=0$.

c) First note that $-\xi_e(x)\le\xi_e(-x)\le0$. Indeed, from subadditivity of $\xi_e$ we have
$0=\xi_e(x+(-x))\le\xi_e(x)+\xi_e(-x)$ which proves the first inequality. Obviously, $-x\notin P$ which
proves the second inequality (except for $x=0$, but this case is trivial). Now, we have
$|\xi_e(-x)|=-\xi_e(-x)\le\xi_e(x)$.

\end{pf}

The main result of this section is the following Theorem.

\begin{Theorem}\label{equinorms} Let $(E,P)$ be an ordered tvs, and let
$e,e'\in\intt P$. Then $\|\cdot\|_e$ and $\|\cdot\|_{e'}$ are mutually equivalent norms. Moreover, it holds
\begin{equation}\label{equi}
\frac1{\xi_e(e')}\|x\|_e\le\|x\|_{e'}\le\xi_{e'}(e)\|x\|_e.
\end{equation}
Also, the constants in (\ref{equi}) are the best possible.
\end{Theorem}

\begin{Remark} If $x\in P$ then by Lemma \ref{l1} c), $||x||_e=\xi_e(x)$ and (\ref{equi}) is reduced to
$$\frac1{\xi_e(e')}\xi_e(x)\le\xi_{e'}(x)\le\xi_{e'}(e)\xi_e(x).$$
\end{Remark}

\begin{pf} The first inequality follows from the second by exchanging roles of $e$ and $e'$. So, we shall
prove only the second.

First, note that $\xi_e(e)=\xi_{e'}(e')=1$ and $\xi_{e'}(e)$, $\xi_e(e')>0$. Let $\lambda\ge0$ and
$\lambda\ge\xi_{e'}(e)\xi_e(x)$. Then we have
$$\lambda e'-x=\frac\lambda{\xi_{e'}(e)}(\xi_{e'}(e)e'-e)+\frac\lambda{\xi_{e'}(e)}e-x.\nonumber
$$
The first summand belongs to $P$, since $\lambda\ge0$, whereas the second belongs to $P$, since
$\lambda/\xi_{e'}(e)\ge\xi_e(x)$. Therefore $\lambda e'-x\in P$. Thus, we got the implication
\begin{equation}\label{lambda>0}
\lambda\ge0\:\wedge\:\lambda\ge\xi_{e'}(e)\xi_e(x)\Longrightarrow\lambda\ge\xi_{e'}(x).
\end{equation}

Similarly, for $\lambda\le0$ and $\lambda\le-\xi_{e'}(e)\xi_e(-x)$, we have
$$-(\lambda
e'-x)=-\frac\lambda{\xi_{e'}(e)}(\xi_{e'}(e)e'-e)+\left(-\frac\lambda{\xi_{e'}(e)}e+x\right).\nonumber
$$
The first summand belongs to $P$, since $\lambda\le0$, whereas the second belongs to $P$, since
$-\lambda/\xi_{e'}(e)\ge\xi_e(-x)$. Thus, $-(\lambda e'-x)\in P$, implying $\lambda e'-x\notin P$; hence
$\lambda\le\xi_{e'}(x)$. We got the second implication:
\begin{equation}\label{lambda<0}
\lambda\le0\:\wedge\:\lambda\le-\xi_{e'}(e)\xi_e(-x)\Longrightarrow\lambda\le\xi_{e'}(x).
\end{equation}

Continuing the proof, we consider three cases $\xi_{e'}(x)>0$, $\xi_{e'}(x)<0$ and $\xi_{e'}(x)=0$.

First case: $\xi_{e'}(x)>0$. If $\xi_e(x)\le0$, then choosing $\lambda=0$ in (\ref{lambda>0}) we get
$\xi_{e'}(x)\le0$ - a contradiction; therefore $\xi_e(x)>0$. Applying, once again (\ref{lambda>0}), now for
$\lambda=\xi_{e'}(e)\xi_e(x)$, we conclude that
\begin{equation}\label{t11}
0<\xi_{e'}(x)\le\xi_{e'}(e)\xi_e(x)\le\xi_{e'}(e)\|x\|_e.
\end{equation}

Second case: $\xi_{e'}(x)\le0$. If $\xi_e(-x)\le0$, then by (\ref{lambda<0}), choosing $\lambda=0$ we obtain
$\xi_{e'}(x)\ge0$ - a contradiction; therefore $\xi_e(-x)>0$. Applying (\ref{lambda<0}) again, for
$\lambda=-\xi_{e'}(e)\xi_e(-x)$, we get
\begin{equation}\label{t12}
0>\xi_{e'}(x)\ge-\xi_{e'}(e)\xi_e(-x)\ge-\xi_{e'}(e)\|x\|_e.
\end{equation}

Combining (\ref{t11}) and (\ref{t12}) we have $|\xi_{e'}(x)|\le\xi_{e'}(e)\|x\|_e$, for $\xi_{e'}(x)\neq0$,
whereas for $\xi_{e'}(x)=0$ this inequality is trivial. Apply the last inequality to $-x$ instead of $x$ to
obtain $\|x\|_{e'}\le\xi_{e'}(e)\|x\|_e$.

To see that constants in (\ref{equi}) are the best possible, choose $x=e$ and $x=e'$.
\end{pf}

\begin{Corollary} Let $(E,P)$ be an ordered tvs. Then:

a) There is a (real) norm on $E$ that generates topology contained in the initial one;

b) The initial topology is Hausdorff;
\end{Corollary}

\begin{pf} a) By Lemma \ref{l1}, $\|\cdot\|_e$ is a norm on $E$. The corresponding topology is generated by balls,
and taking into account that $x\mapsto a+x$ and $x\mapsto\alpha x$ ($a\in E$, $\alpha\in\mathbf R$) are
continuous functions, it is enough to prove that $\{x\in E~|~\|x\|_e<1\}$ is open in the initial topology.
However it immediately follows, since the latter set is equal to $\xi_e^{-1}(-1,1)\cap-\xi_e^{-1}(-1,1)$, and
$\xi_e$ is continuous.

b) Follows from the previous part.
\end{pf}

\begin{Corollary} Let $E$ be a topological vector space. If $E$ does not contain nontrivial
(i.e.~$\neq\emptyset,E$) convex open sets then there is no solid cone in $E$. For instance, $L^p(0,1)$ for
$0<p<1$ has not any solid cone.
\end{Corollary}

\begin{pf} By previous Corollary, the unit ball in the norm $\|\cdot\|_e$ is open. It is also (obviously)
convex, so the result follows. The space $L^p(0,1)$, $0<p<1$ does not contain any nontrivial open convex set
by \cite[\S1.47]{Rudin}.
\end{pf}

\section{TVS convergence}

\begin{Definition} Let $X$ be a nonempty set, and let $(E,P)$ be an ordered tvs.
If the function $d:X\times X\to P$ satisfies the following conditions:

\begin{enumerate}
\item $d(x,y)=0$ if and only if $x=y$;

\item $d(y,x)=d(x,y)$, for all $x,y\in X$;

\item $d(x,z)\le d(x,y)+d(y,z)$ (in the sense of Definition \ref{poredak}), for all $x,y,z\in X$,
\end{enumerate}
then the pair $(X,d)$ (or quadruple $(X,E,P,d)$ if we want to emphasize ordered tvs) is called TVS-cone
metric space.
\end{Definition}

The following definition is well known \cite{HZ}, \cite{DU}

\begin{Definition} Let $(X,d)$ be a TVS-cone metric space, let $(x_n)\subseteq X$ be a sequence, and let
$x\in X$.

\begin{enumerate}
\item We say that $(x_n)$ TVS converges to $x$ if for any $c\in\intt P$ there is a positive integer
$n_0$, such that for all $n\ge n_0$ we have $d(x_n,x)\ll c$;

\item We say that $(x_n)$ is a TVS Cauchy sequence, if for all $x\in\intt P$ there is a positive
integer $n_0$ such that for all $m,n\ge n_0$ we have $d(x_m,x_n)\ll c$;

\item We say that $(X,d)$ is TVS-cone complete if every TVS Cauchy sequence is TVS convergent.
\end{enumerate}
\end{Definition}

\begin{Remark} It is easy to see, using Lemma \ref{basic} - 2) that "for all $c\in\intt P$" can be replaced
by "for all $c=\e e$, $\e$ runs through $\mathbf R^+$ and $e\in\intt P$ is fixed". Thus, TVS convergence is
controlled by $\mathbf R$ instead of $\intt P$.
\end{Remark}

In a similar manner we can define TVS continuous and TVS uniformly continuous functions.

\begin{Definition} Let $(X_1,d_1)$ and $(X_2,d_2)$ be TVS-cone metric spaces over an ordered tvs $(E,P)$.
\begin{enumerate}
\item We say that $f:X_1\to X_2$ is TVS continuous at $a\in X_1$, if for all $e\in\intt P$ there exists $c\in\intt
P$ such that for all $x\in X_1$
$$d_1(a,x)\ll c\mbox{ implies }d_2(f(a),f(x))\ll e.$$
We say that $f$ is TVS continuous if it is TVS continuous at all $a\in X_1$;

\item Similarly, we say that $f$ is TVS uniformly continuous if for all $e\in\intt P$ there exists $c\in\intt P$
such that for all $x,y\in X_1$
$$d_1(x,y)\ll c\mbox{ implies }d_2(f(x),f(y))\ll e.$$
\end{enumerate}
\end{Definition}

The next Proposition 1)-4) is a refinement of a result proved in \cite[Lemma 2]{AK} and \cite[Theorems 2.1
and 2.2]{DU} which have additional assumption that $E$ is locally convex.

\begin{Proposition}\label{Cauchy-equi}
Let $(X,d)$ and $(X',d')$ be TVS-cone metric spaces over an ordered tvs $(E,P)$. Then:

1) The function $d_e:X\times X\rightarrow [0,+\infty)$ defined by $d_e=\xi_e\circ d$ is a metric;

2) The sequence $(x_n)$ TVS-cone converges to $x$, if and only if $d_e (x_n,x)\to0$ for all $e\in\intt P$ and
this holds if and only if $d_e(x_n,x)\to 0$ for some $e\in\intt P$;

3) The sequence $(x_n)$ is a TVS-cone Cauchy sequences if and only if $(x_n)$ is a Cauchy sequences (in usual
sense) in $(X,d_e)$ for all $e\in\intt P$ and this holds if and only if $(x_n)$ is a Cauchy sequence in
$(X,d_e)$ for some $e\in\intt P$;

4) $(X,d)$ is a TVS-cone complete metric space, if and only if all $(X,d_e)$ is complete metric space and
this holds if and only if $(X,d_e)$ is Cauchy complete for some $e\in\intt P$.

5) $f:X\to X'$ is TVS continuous if and only if it is continuous in all $d_e$, $d'_e$, $e\in\intt P$;

6) $f:X\to X'$ is TVS uniformly continuous if and only if it is uniformly continuous in all $d_e$, $d'_e$,
$e\in\intt P$.
\end{Proposition}

\begin{pf}

1) By Lemma \ref{l1}, we have $d_e(x,y)=\|d(x,y)\|$ and desired properties follows immediately from
properties of $\|\cdot\|_e$ and $d$.

By Theorem \ref{equinorms} and Remark after it, all metrics $d_e$ are mutually equivalent. Thus, we can
assume in the rest of the proof that the equivalence between "in all $d_e$" and "in some $d_e$" has been
already proved.

2) Let $x_n$ TVS converges to $x$, and let $e\in\intt P$, $\e>0$ be arbitrary. Then $\e e\in\intt P$. So
$d(x_n,x)\ll\e e$ for $n$ large enough. By Lemma \ref{sf} we obtain $\xi_e(d(x_n,x))\le\e$ for same $n$.

Let us suppose that $x_n\to x$ in all metrics $d_e$, let $e\in\intt P$, and let $\e<1$. Then
$d_e(x_n,x)=\xi_e(d(x_n,x))\to0$, and therefore there exists a positive integer $n_0$ such that for all $n\ge
n_0$ we have
$$\inf\{r\:|\:d(x_n,x)\in re-P\}=\xi_e(d(x_n,x))\le\varepsilon/2,$$
which implies that there exists $r<\varepsilon$ such that $d(x_n,x)\in re-P$, i.e.~$re-d(x_n,x)\in P$,
i.e.~$d(x_n,x)\le re\ll e$.

3) The proof is the same as the proof of the previous part of the statement. Just replace $d(x_n,x)$ by
$d(x_n,x_m)$.

4) It follows from conclusions 2) and 3).

5) Let $f:X\to X'$ be TVS continuous, and let $e\in\intt P$ be fixed. Given $\e>0$, we have $\e e\in\intt P$;
hence there is a $c\in\intt P$ such that $d(a,x)\ll c$ implies $d'(f(a),f(x))\ll\e e$. Choose $\delta>0$ such
that $\delta e\ll c$ (Lemma \ref{basic} - 2)). Then $\xi_e(d(a,x))=d_e(a,x)<\delta$ implies $d(a,x)\le\delta
e\ll c$, and therefore $d'(f(a),f(x))\ll\e e$ which leads to $\xi_e(d'(f(a),f(x)))\le\e$.

Conversely, let $e\in\intt P$ be arbitrary. We shall use the fact that $f$ is continuous in $d_e$. For
$\e=1/2$, there is $\delta>0$ such that for all $x$, $\xi_e(d(a,x))\le\delta$ implies
$\xi_e(d'(f(a),f(x)))<1/2$. Choose $c=\delta e$; if $d(a,x)\ll c=\delta e$ then $d_e(a,x)\le\delta$ which
implies $d'_e(f(a),f(x))=\xi_e(d'(f(a),f(x)))<1/2$ which leads to $d'(f(a),f(x))\le e/2\ll e$.

6) The proof is the same as that of 5).
\end{pf}

\begin{Proposition}\label{onemetric}
Let $(X,d)$ be a TVS-cone metric space (over an ordered tvs $(E,P)$). Then there exists a metric $d^*$ on $X$
such that notions "being convergent", "being Cauchy sequence" are equivalent in $d$ and $d^*$.
\end{Proposition}

\begin{pf} Pick some $e\in\intt P$ and consider $d^*=d_e$. If $(x_n)$ is TVS-convergent to $x$
(resp.~TVS-Cauchy sequence) then $(x_n)$ is convergent to $x$ (resp.~is a Cauchy sequence) by previous
Theorem. This proves one direction.

For the other, let $(x_n)$ converges to $x$ (resp.~be a Cauchy sequence) in the metric $d_e$, and let
$e'\in\intt P$ be arbitrary. By Lemma \ref{l1} part c), $d_e(x,y)=\xi_e(d(x,y))=\|d(x,y)\|_e$ and
$d_{e'}=\|d(x,y)\|_{e'}$. By Theorem \ref{equinorms}, norms $\|\cdot\|_e$ and $\|\cdot\|_{e'}$ are
equivalent. So, $(x_n)$ converges to $x$ (resp.~is a Cauchy sequence) in all metrics $d_{e'}$, $e'\in\intt
P$; hence so it is in TVS sense.
\end{pf}

\begin{Remark} The previous Proposition was proved in \cite{KP} in the special case of Banach spaces using
quite different methods, in \cite{AR} for locally convex spaces also using different methods. It was also
proved in \cite{KRR1} using Minkowsky functional, under additional assumption that the set $[-e,e]=\{x\in
E~|~-e\le x\le e\}$ is absolutely convex. However, the authors in \cite{KRR1} miss to specify whether they
work with complex or real scalars. If scalars are real (or if they are complex, but the tvs is considered as
a real linear space) then $[-e,e]$ is absolutely convex, indeed. It follows, for instance, from Lemma
\ref{l1}. In this case Minkowsky functional $q_e$ in the proof of \cite[Theorem3.2]{KRR} is equal to the norm
$||\cdot||_e$. However, \cite[Remark to Lemma 1]{AK} shows that $[-e,e]$ might not be absolutely convex, if
$E$ is considered as a linear space over complex field. In more details, $E=\mathbf C$ with standard
topology, $P=\{z\in\mathbf C~|~\im z\ge2|\re z|\}$, $e=i$, $x=\delta i$ for some $1/2<\delta<1$. Then
$x\in[-e,e]$, but $ix\notin[-e,e]$. Hence $[-e,e]$ is not balanced and therefore it is not absolutely convex.
\end{Remark}

The previous Proposition establishes that TVS convergence are equivalent to the convergence in some (usual)
metric space. However, in many applications such as fixed point results, one uses explicitly the metric. The
next Proposition establishes the connection with inequalities in TVS metric and metrics $d_e$.

\begin{Proposition}\label{nejequi}
Let $(E,P)$ be an ordered tvs, and let $x,y\in P$ Then $x\le_P y$ if and only if for all $e\in\intt P$ there
holds
\begin{equation}\label{ineq}
\xi_e(x)\le\xi_e(y).
\end{equation}
\end{Proposition}

\begin{pf} One part follows from the monotonicity of $\xi_e$.

For the other, consider first the case $y\in\intt P$. Then apply (\ref{ineq}) for $y=e$ to obtain
$\xi_y(x)\le\xi_y(y)=1$; hence $1\in M_{y,x}$ implying $y-x\in P$. If $y\in\partial P$ then we can prove the
statement by limit argument. Let $z\in\intt P$ be arbitrary. Obviously, $\intt P\ni z/n\to0$. Choosing
$e=y+z/n$ (which $\in\intt P$ by Lemma \ref{basic} -1)) in (\ref{ineq}) we have
$\xi_{y+z/n}(x)\le\xi_{y+z/n}(y)$, or equivalently
\begin{equation}\label{ineq1}\xi_{y+z/n}(y+z/n)-x\in P.
\end{equation}
Thus, the first summand in the next formula belongs to $P$.
\begin{equation}\label{ineq2}
y+z/n-x=\left(\xi_{y+z/n}(y)(y+z/n)-x\right)+(1-\xi_{y+z/n}(y))(y+z/n).
\end{equation}
On the other hand, since $1\cdot (y+z/n)-y=z/n\in P$, it follows $1\in M_{y+z/n,y}$ implying
$\xi_{y+z/n}(y)\le1$. Therefore the second summand in (\ref{ineq2}) also belongs to $P$. So, $y+z/n-x\in P$
and the result follows since $P$ is closed.
\end{pf}

The results from this section can be summarized from the categorial point of view in the following way.

Given a fixed ordered tvs $(E,P)$, consider the categories $\mathcal E_1$, $\mathcal E_2$ and $\mathcal E_3$
whose objects are TVS-cone metric spaces and whose morphisms are, respectively, TVS continuous mappings, TVS
uniformly continuous mappings and mappings that not increase the distance (i.e, the mappings $f:X\to X'$ such
that $d'(f(x),f(y))\le_P d(x,y)$). On the other hands, consider the categories $\mathcal M_1$, $\mathcal M_2$
and $\mathcal M_3$ whose objects are pairs $(X,D)$, where $X$ ia a nonempty set and $D=\{d_e~|~e\in\intt P\}$
is a family of mutually equivalent norms, and whose morphisms are, respectively, mappings that are continuous
in all $d_e$, uniformly continuous in all $d_e$, and distance not increasing with respect to all $d_e$.

\begin{Corollary}\label{categoty} We have:

1) The mapping $F_i:\mathcal E_i\to\mathcal M_i$, $i=1,2,3$, that sends object $(X,d)$ to
$(X,\{d_e~|~e\in\intt P\})$ and morphism $f:X_1\to X_2$ to the same $f$ is a covariant functor;

2) $\mathcal E_i$ is a proper subcategory of $\mathcal M_i$.
\end{Corollary}

\begin{pf} 1) is a direct consequence of Propositions \ref{Cauchy-equi}, \ref{onemetric} and \ref{nejequi}
and Theorem \ref{equinorms};

2) If $f:X_1\to X_2$ is an isomorphism in $\mathcal M_i$, then $f$ is bijective, and $f^{-1}$ has the same
properties as $f$. For $i=1$, $f$, $f^{-1}$ are continuous in all $d_e$ and hence in TVS metric $d$ by
Proposition \ref{Cauchy-equi} - 5). For $i=2$, $f$, $f^{-1}$ are uniformly continuous in all $d_e$ and hence
in TVS metric $d$ by Proposition \ref{Cauchy-equi} - 6). Finally, for $i=3$, we have $d_e(fx,fy)\le d_e(x,y)$
and so for $f^{-1}$; hence by Proposition \ref{nejequi} we have $d(fx,fy)\le_P d(x,y)$, and also
$d(f^{-1}x,f^{-1}y)\le_P d(x,y)$. By antisymmetry we have that $f$ is a TVS isometry.

Moreover, there are some objects in $\mathcal M_i$ that can not be obtained as an image of the functor $F_i$.
For instance, the objects that do not satisfy $d_{\lambda e}=(1/\lambda)d_e$ for some $e\in\intt P$.
\end{pf}

As a usual relationship between more and less general theory, we have

\begin{Corollary} There is no generalization of results from metric spaces (equipped with a family of
mutually equivalent norms) to TVS-cone metric spaces. Any statement concerning TVS-cone metric spaces is
either special case of an metric spaces result or a new theorem that can not proved in metric spaces.
\end{Corollary}

\begin{Remark} Result analogous to the Corollary \ref{categoty} can be obtained for other "metric-like" structures such
as semi-metric spaces etc. We consider there is no need to list all of them.
\end{Remark}

\begin{Remark} Despite of the previous two Corollaries there are some nontrivial questions. For instance, it
is clear that it must be $d_{\lambda e}=(1/\lambda)d_e$, but it is nontrivial question how to express
$d_{e_1+e_2}$ in terms of $d_{e_1}$ and $d_{e_2}$, since this extremely depends on the "shape" of the cone
$P$. Also, nontrivial question is can the space $E$ together with cone $P$ be recovered by the family of
metrics $d_e$.
\end{Remark}

\section{Extensions of fixed point results}

First we give the statement that allows automatical derivations of fixed point results from (usual) metric
spaces to TVS-cone metric spaces, as their special cases. Its proof is direct consequence of Corollary
\ref{categoty}.

\begin{Theorem}\label{trivgen} Let $(E,P)$ be an ordered tvs, let $(X,d)$ be TVS-cone metric space, and let $f:X\to X$ be a
mapping. The statement
\begin{quote}If $f$ satisfies conditions $A_i\le_P B_i$, $i=1,2,\dots,n$, where $A_i$,
$B_i$ are terms made by constants, variables (from $X$, $\mathrm R$), multiplication $\mathbf R\times E\to E$
and $d$, then $f$ has a unique fixed point\end{quote}
is equivalent to the statement
\begin{quote}If $f$ satisfies conditions $A^e_i\le B^e_i$, $i=1,2,\dots,n$ for all $e\in\intt P$ where terms $A^e_i$,
$B^e_i$ are obtained by formal substitution of $d$ by $d_e$, then $f$ has a unique fixed point.
\end{quote}
\end{Theorem}

\begin{Remark} Thus, many generalizations of fixed point results from (usual) metric spaces to TVS-cone
metric spaces are not really generalizations. On the contrary they are weaker results, since they suppose
that contractivity condition is satisfied for a family of metrics instead of a single one.
\end{Remark}

However, if contractivity conditions contain constants or variables from $P$ then automatic extension using
Theorem \ref{trivgen} is not possible. So we need to reformulate such conditions. To illustrate this, we
prove a new common fixed point result for metric spaces, that is a generalizations of K.~M.~Das and
K.~V.~Naik \cite{DN} common fixed point theorem, and show how it can be reduced to TVS-cone metric spaces in
few steps.

First, recall some standard terminology and notations from fixed point theory.

Let $X$ be a nonempty set and let $f:X\rightarrow X$ be an arbitrary mapping. The element $x\in X$ is a fixed
point for $f$ if $x=f(x)$.

Let $X,Y$ be a nonempty sets, $f,g:X\rightarrow Y$ and $f(X)\subseteq g(X)$. Choose a point $x_1\in X$ such
that $f(x_0)=g(x_1)$. Containing this process, having chosen $x_n\in X$, we obtain $x_{n+1}\in X$ such that
$f(x_n)=g(x_{n+1})$. $f(x_n)$ is called Jungck sequence with initial point $x_0$. Note that Jungck sequence
might not be determined by its initial point $x_0$.

Let $X,Y$ be a nonempty sets and $f,g:X\rightarrow Y$. If $f(x)=g(x)=y$ for some $y\in Y$ then $x$ is called
a coincidence point of $f$ and $g$, and $y$ is called a point of coincidence of $f$ and $g$.

Let $X$ be a nonempty set and $f,g:X\rightarrow X$. $f$ and $g$ is weakly compatible self mappings if they
commute at their coincidence point.

\begin{Lemma} (\cite{AJ}) Let $X$ be a nonempty set and let $f,g:X\rightarrow X$ be weakly compatible self mappings. If
$f$ and $g$ have the unique point of coincidence $y=f(x)=g(x)$, then $y$ is the unique common fixed point of
$f$ and $g$.
\end{Lemma}

We shall consider certain classes of function described in the next Definition.

\begin{Definition}
By $\Phi $ we denote the set of all real functions $\varphi:[0,\infty) \rightarrow [0,\infty)$ which have the
following properties:

(a) $\varphi (0)=0$;

(b) $\varphi(x) < x$ for all $x>0$;

(c) $\lim_{x\rightarrow\infty} (x-\varphi(x))=\infty$.

\noindent Define $\Phi_1=\{\varphi\in\Phi:\varphi \mbox{ is monotone nondecreasing and } \overline{\lim}_{t
\rightarrow r+} \varphi(t)< r\mbox{ for any } r>0\}$ and $\Phi_2=\{\varphi\in\Phi:\overline{\lim}_{t
\rightarrow r} \varphi(t)< r\mbox{ for any } r>0\}$.\end{Definition}

Two following two Lemmas were proved in \cite{ARK}.

\begin{Lemma} \label{ark1} Let $\varphi\in\Phi_2$. Then there exists
$\psi\in\Phi_1$ such that $$\varphi (x)\le\psi (x) <x,$$ for each $x>0$. \end{Lemma}

\begin{Lemma} \label{ark2} Let $\varphi_1,\ldots,\varphi_n \in\Phi_1$. Then there exists
$\psi\in\Phi_1$ such that $$\varphi_k (x)\le\psi (x) <x,$$ for each $1\le k\le n$ and $x>0$. \end{Lemma}

Before stating the result, we make a convention to abbreviate $f(x)$ and $g(x)$ in order to avoid too much
parenthesis.

\begin{Theorem}\label{mainfixed}
Let $X$ be a nonempty set, let $(Y,d)$ be a metric space and let $f,g:X\rightarrow Y$ be two mappings.
Suppose that the range of $g$ contains the range of $f$ and that $g(X)$ is a complete subspace of $Y$. If
there exists $\varphi_1,\varphi_2,\varphi_3,\varphi_4,\varphi_5\in\Phi_2$ such that
\begin{eqnarray}\label{contr} d(fx,fy) &\le &\max\{\varphi_1(d(gx,gy)),\varphi_2(d(gx,fx)),
\varphi_3(d(gy,fy)),\nonumber \\
&& \varphi_4(d(gx,fy)),\varphi_5(d(fx,gy)\},
\end{eqnarray}
for any $x,y\in X$, then there exists $z\in Y$ which is the limit of every Jungck sequence defined by $f$ and
$g$. Further, $z$ is the unique point of coincidence of $f$ and $g$. Moreover, if $X=Y$ and $f$, $g$ are
weakly compatible then $z$ is the unique common fixed point for $f$ and $g$.

\end{Theorem}

\begin{pf} We shall, first, reduce the statement to the case $\varphi_1=\dots=\varphi_5\in\Phi_1$. Indeed,
from Lemma \ref{ark1} it follows that there exist functions $\varphi_k^*\in\Phi_1$ such that
$\varphi_k^*\in\Phi_1$ such that
$$\varphi_k(x)\le\varphi_k^* (x) <x,$$
for each $x>0$($1\le k\le 5$), whereas, from Lemma \ref{ark2} it follows that there exists a real function
$\varphi\in\Phi_1$ such that:
$$\varphi_k^*(x)\le\varphi (x) <x,~(1\le k\le 5)\mbox{ for
each }x>0,$$
which implies
\begin{eqnarray}d(fx,fy) &\le
&\max\{\varphi(d(gx,gy)),\varphi(d(gx,fx)),
\varphi(d(gy,fy)),\nonumber \\
&& \varphi(d(gx,fy)),\varphi(d(fx,gy))\}.\nonumber
\end{eqnarray}

Thus, we can assume that $\varphi_j=\varphi\in\Phi_1$ for $1\le j\le 5$.

Let $x_0\in X$ be arbitrary and let $(x_n)$ be an arbitrary sequence such that $f(x_n)$ is Jungck sequence
with initial point $x_0$.

Let $d_0=d(fx_0,gx_0)$. We will prove that there exists a real number $r_0>0$ such that:
\begin{equation} r_0-\varphi(r_0)\le d_0\quad\mbox{and}\quad r-\varphi(r)>d_0\quad\mbox{for}\quad
r>r_0.
\end{equation}
Consider the set $D=\{r~|~\forall s>r,\:s-\varphi(s)>d_0\}$ which is nonempty, since $r-\varphi(r)\to\infty$
as $r\to\infty$. Also, it holds $s\in D$, $t>s$ implies $t\in D$, and hence, $D$ is an unbounded interval.
Set $r_0=\inf D$. For each positive integer $n$, there is $r_n\notin D$ such that $r_0-1/n<r_n$, and
therefore, there is $r_0\geq s_n>r_n>r_0-1/n$ such that $s_n-\varphi(s_n)\leq d_0$. Since $\varphi$ is
nondecreasing, we have $\varphi(s_n)\leq\varphi(r_0)$, implying $s_n-\varphi(r_0)\leq d_0$. Taking a limit as
$n\to\infty$ we get $r_0-\varphi(r_0)\leq d_0$.

For any $j\ge0$, define $\mathcal{O}_n(x_j)=\{f(x_k)~|~k=j,j+1,j+2,\ldots,j+n\}$,
$\mathcal{O}(x_j)=\{f(x_k)~|~k=j,j+1,j+2,\ldots\}$. Also, let $\delta(A)$ denote the diameter of $A$.

Next, we prove that for all positive integers $k,n$ there holds
\begin{equation}\label{Odiam}
\delta(\mathcal O_n(x_k))\leq\varphi(\delta(\mathcal O_{n+1}(x_{k-1}))).
\end{equation}

Since $\varphi$ is nondecreasing, it commutes with $\max$, and for $k\le i,j\le k+n$, we have
\begin{eqnarray}  d(fx_i,fx_j) \le \varphi(\max\{d(gx_i,gx_j),d(gx_i,fx_i),
d(gx_j,fx_j),\nonumber\\d(gx_i,fx_j),d(gx_j,fx_i)\})
=\nonumber \\
=\varphi(\max\{d(fx_{i-1},fx_{j-1}),d(fx_{i-1},fx_i),
d(fx_{j-1},fx_j),\nonumber\\d(fx_{i-1},fx_j),d(fx_{j-1},fx_i)\})\le\nonumber\\
\varphi(\delta(\mathcal{O}_{n+1}(x_{k-1}))).\nonumber
\end{eqnarray}

By induction, from (\ref{Odiam}) we obtain
\begin{equation}\label{Odiam2}
\delta(\mathcal O_n(x_k))\le\varphi^l(\delta(\mathcal O_{n+l}(x_{k-l}))).
\end{equation}

For $1\le i,j\le n$ we have $fx_i,fx_j\in\mathcal O_{n-1}(x_1)$, and hence, by (\ref{Odiam})
$$d(fx_i,fx_j)\le\delta(\mathcal O_{n-1}(x_1))\leq\varphi(\delta(\mathcal O_n(x_0)))<\delta(\mathcal
O_n(x_0)).$$
Therefore, there is $1\le k\le n$, such that
\begin{eqnarray}\delta(\mathcal O_n(x_0))=d(fx_0,fx_k)\leq d(fx_0,fx_1)+d(fx_1,fx_k)\leq\nonumber
\\\leq d_0+\delta(\mathcal
O_{n-1}(x_1))\leq d_0+\varphi(\delta(\mathcal O_n(x_0)).\nonumber
\end{eqnarray}
Hence we get
$$\delta (\mathcal{
O}_n(x_0))-\varphi(\delta (\mathcal{O}_n(x_0)))\le d_0$$
which implies $\delta (\mathcal{O}_n(x_0))\le r_0$, implying
\begin{equation}\label{bounded}
\delta(\mathcal O(x_0))=\sup_n\delta(\mathcal O_n(x_0))\le r_0.
\end{equation}
Hence all Jungck sequences defined by $f$ and $g$ are bounded.

Now, we shall prove that our Jungck sequence is a Cauchy sequence. Let $m>n$ be positive integers. Then
$fx_n,fx_m\in\mathcal O_{m-n+1}(x_n)$, and using (\ref{Odiam2}) (with $l=n$) and (\ref{bounded}) we get
$$d(fx_n,fx_m)\le\delta(\mathcal O_{m-n+1}(x_n))\le\varphi^n(\delta(\mathcal O_{m+1}(x_0)))\le\varphi^n(r_0)\to0,$$
as $m,n\to\infty$. Since $f(X)\subseteq g(X)$, and $g(X)$ is complete, it follows that $fx_n$ is convergent.
Let $y\in X$ be its limit.

Clearly $y\in g(X)$, so there is $z\in X$ such that $g(z)=y$. Let us prove that $f(z)$ is also equal to $y$.
By (\ref{contr}) we have
\begin{eqnarray}d(fx_n,fz)\le\varphi(\max\{d(gx_n,gz),d(gx_n,fx_n),d(gz,fz),\nonumber\\
d(gx_n,fz),d(gz,fx_n)\})=\nonumber\\
=\varphi(\max\{d(fx_{n-1},y),d(fx_{n-1},fx_n),d(y,fz),\nonumber\\d(fx_{n-1},fz),d(y,fx_n)\}).\nonumber
\end{eqnarray}
If $n\to\infty$, then the left hand side in the previous inequality tends to $d(y,fz)$, the first, the second
and the fifth argument of $\max$ tend to $d(y,y)=0$, whereas the third and the fourth tend to $d(y,fz)$.
Thus, we have
$$d(y,fz)\le\varphi(d(y,fz)),$$
which is impossible, unless $d(y,fz)=0$.

Finally, we prove that the point of coincidence is unique. Suppose that there is two points of coincidence
$y$ and $y'$ obtained by $z$ and $z'$, i.e.~$fz=gz=y$ and $fz'=gz'=y'$. Then by (\ref{contr}) we have
\begin{eqnarray}d(y,y')=d(fz,fz')\le\varphi(\max\{d(gz,gz'),d(gz,fz),d(gz',fz'),\nonumber\\
d(gz,fz'),d(gz',fz)\})=\nonumber\\
\varphi(\max\{d(y,y'),0,0,d(y,y'),d(y',y)\})=\varphi(d(y,y'))<d(y,y'),\nonumber
\end{eqnarray}
unless $d(y,y')=0$. Since every Jungck sequence converges to some point of coincidence, and the point of
coincidence is unique, it follows that all Jungck sequences converge to the same limit.

Let $X=Y$ and let $f$, $g$ be weakly compatible. By Lemma 3 we get that $y=z$ which is the unique common
fixed point of $f$ and $g$.
\end{pf}

\begin{Remark} A special case of the previous statement where $\varphi_i\in\Phi_1\cup\Phi_2\subseteq\Phi_2$, $X=Y$ and
$g(x)\equiv x$ was proved in \cite{ARK}. It generalizes earlier results \cite{I}, \cite{D}, \cite{C2}.
\end{Remark}

In the sequel we want to reformulate this result to TVS-cone metric spaces. We need counterparts of classes
$\Phi$, $\Phi_2$.

Let $(E,P)$ be an ordered tvs. For a function $g:[t_0,+\infty)\to P$, we say that
$\lim\limits_{t\to+\infty}g(t)=\infty$, if there exists $c\in\intt P$ and a scalar function
$h:[t_0,+\infty)\to\mathbf R$, such that
\begin{equation}\label{infty}
\lim\limits_{t\to+\infty}h(t)=+\infty\quad\mbox{and}\quad h(t)c\ll g(t).
\end{equation}
Note that this is equivalent to the property that for all $c\in\intt P$ there exists $M>0$ such that $t>M$
implies $g(t)\gg c$, as well as to the statement that for all $c\in\intt P$ there holds (\ref{infty}) for
some $h$ (Lemma \ref{basic}).

By $\Psi$ we denote the set of all functions $\psi:P\rightarrow P$ which have the following properties:

(a) $ \psi(0)=0$;

(b) $(I-\psi)(\intt P)\subseteq \intt P$;

(c)  $\lim_{t\rightarrow\infty} tx-\psi (tx)=\infty$, for any $x\in P\setminus \{0\}$.

For a function $\psi\in\Psi$ we say that $\psi\in\Psi_2$ if for any $x\in \mathrm{int}P$ and for each
$(x_n)\subseteq \mathrm{int}P$ such that $x_n\rightarrow x$, there exists a positive integer $n_0$ such that
$n>n_0$ implies $\psi(x_n)\le (1-\varepsilon)x$, where $\varepsilon$ does not depend on the choice of the
sequence $(x_n)$.

Next result is a generalization of earlier results presented in \cite{AK},\cite{IR1}, \cite{KRR}, \cite{KRR1}
and \cite{KRR2} (some of them also can be obtained "automatically" by Theorem \ref{trivgen}).

\begin{Theorem}\label{last}
Let $X$ be a nonempty set, $(Y,d)$ be a TVS-cone metric space and $f,g:X\rightarrow Y$ two mappings. Suppose
that the range of $g$ contains the range of $f$ and $g(X)$ is a complete subspace of $Y$.  Assume that there
exists $\psi_1,\psi_2,\psi_3,\psi_4,\psi_5\in\Psi_2$ such that for any $x,y\in X$ there exists
\begin{eqnarray}\label{prva1} u\in \{\psi_1(d(g(x),g(y))),\psi_2(d(g(x),f(x))),
\psi_3(d(g(y),f(y))),\nonumber \\
\psi_4(d(g(x),f(y))),\psi_5(d(f(x),y))\}\end{eqnarray}
such that
\begin{equation}\label{druga}d(f(x),f(y)) \le u,\end{equation}
then there exists $z\in Y$ which is the limit of every Jungck sequence defined by $f$ and $g$. Further $z$ is
the unique point of coincidence for $f$ and $g$. Moreover, if $X=Y$ and $f$, $g$ are weakly compatible  then
$y$ is unique common  fixed point for $f$ and $g$.
\end{Theorem}

\begin{pf} Let $e\in \mathrm{int} P$ be any. For $k\in \{1,\ldots, 5\}$ define
$\varphi_k:[0,\infty) \rightarrow [0,\infty)$ by
$$\varphi_k(t)=\xi_e(\psi_k(te))=\|\psi_k(te)\|_e.$$
Then

(a) $\varphi_k(0)=0$ since $\xi_e(0)=0$ and $\psi_k(0)=0$.

(b) By Lemma \ref{l1}, taking into account $\psi_k(te)\ll te$, we have $\varphi_k(t)=\|\psi_k(te)\|<t$.

(c) $\lim_{t\to\infty}(t-\varphi_k(t))=\infty$. Indeed, using monotonicity of $\xi_e$ and definition of
tending to $\infty$ in TVS sense, we have $t-\varphi_k(t)=\xi_e(te-\psi_k(te))\ge\xi_e(h(t)c)=h(t)\xi_e(c)$,
for some $h(t)$ that tends to $\infty$ as $t\to\infty$.

Therefore, $\varphi_k\in\Phi$.

(d) If $\psi_k\in\Psi_2$ then $\varphi_k\in\Phi_2$.

Let $t\in (0,+\infty)$ and $t_n\rightarrow t$. From $\psi_k\in\Psi_2$ and $t_ne\rightarrow te$ it follows
that there exists a positive integer $n_0$ such that $n>n_0$ implies $\psi_k(t_ne)\le(1-\varepsilon)
t_ne\ll(1-\varepsilon/2)t_ne$. By Lemma \ref{l1}, it follows that
$\varphi_k(t_n)=\|\psi_k(t_n)\|<(1-\varepsilon/2)t_n$, for $n>n_0$. Hence $\overline{\lim}_{s \rightarrow t}
\varphi(s)\le(1-\varepsilon/2)t< t$ because $t_n$ and $t$ are arbitrary, which implies that
$\varphi_k\in\Phi_2$.

Now, the statement is a special case of Theorem \ref{mainfixed}, since $\varphi_1,\dots,\varphi_5\in\Phi_2$,
and they obey (\ref{contr}). Moreover, for all $e$ we can find the corresponding
$\varphi_1^e,\dots,\varphi_5^e$.
\end{pf}

\begin{Definition} If $A:E\rightarrow E$ is a one to one function such that
$A(P)=P$, $(I-A)$ is one to one and $(I-A)(P)=P$ then we say that $A$ is contractive operator.
\end{Definition}

Next corollary include results presented in  \cite{IR1}, \cite{KRR}, \cite{C2}, \cite{HZ} and \cite{RH}.

\begin{Corollary} Let $X$ be a nonempty set, let $(Y,d)$ be a TVS-cone metric space and let
$f,g:X\rightarrow Y$ be two mappings. Suppose that the range of $g$ contains the range of $f$ and $g(X)$ is a
complete subspace of $Y$.  Assume that there exists contractive bounded linear operators
$A_1,A_2,A_3,A_4,A_5$ such that for any $x,y\in X$ there exists
\begin{eqnarray} u\in \{A_1(d(g(x),g(y))),A_2(d(g(x),f(x))),
A_3(d(g(y),f(y))),\nonumber \\
A_4(d(g(x),f(y))),A_5(d(f(x),g(y)))\}\nonumber\end{eqnarray}
such that $d(f(x),f(y)) \le u$, then there exists $z\in Y$ which is limit of every Jungck sequence defined by
$f$ and $g$. Further $z$ is the unique point of coincidence for $f$ and $g$. Moreover if $X=Y$ and $f$, $g$
are weakly compatible  then $z$ is the unique common fixed point for $f$ and $g$.
\end{Corollary}

\begin{pf} Let $k\in \{1,\ldots, 5\}$. Then $A_k(0)=0$ since $A_k$ is linear.
From $(I-A_k)(P)=P$ by Open mapping theorem (see \cite{GK}) it follows that $(I-A_k)(\intt P)\subseteq \intt
P$. So $A_k(x)\ll x$ for any $x\in \intt P$. Let $x\in \intt P$, $(x_n)\subseteq \intt P$ and $x_n\rightarrow
x$. Then $\lim A_k(x_n)\ll x$ since $A_k$ is continuous. Also, $tx-A_k(tx)=(x-A_k(x))t$, for any $x \in
P\setminus \{0\}$. Hence $A_k\in \Psi_2$. The statement follows from Theorem \ref{l1}.
\end{pf}

\end{document}